\documentclass[12pt]{article} 
\usepackage{latexsym}
\usepackage{amssymb, amsmath}
\usepackage[left=2.7cm,right=2.7cm,bottom=4cm ]{geometry}

\usepackage{latexsym}
\usepackage{amssymb}
\usepackage{euscript}
\usepackage[dvips]{graphics}
\usepackage{epsf}

\def\mcc{M\raise.5ex\hbox{c}C}
\def\mccarthy{M\raise.5ex\hbox{c}Carthy}

\def\De{\Delta}

\def\={\ = \ }

\def\C{\mathbb C}
\def\R{\mathbb R}

\def\NN{\mathbb N}

\def\dis{\displaystyle}

\def\be{\setcounter{equation}{\value{theorem}} \begin{equation}}
\def\ee{\end{equation} \addtocounter{theorem}{1}}
\def\beq{\begin{eqnarray*}}
\def\eeq{\end{eqnarray*}}

\def\bp{{\sc Proof: }}
\def\ep{{}{\hfill $\Box$} \vskip 5pt \par}

\def\bl{\begin{lemma}}
\def\el{\end{lemma}}
\def\bt{\begin{theorem}}
\def\et{\end{theorem}}
\def\bprop{\begin{prop}}
\def\eprop{\end{prop}}
\def\bd{\begin{definition}}
\def\ed{\end{definition}}
\def\br{\begin{remark}}
\def\er{\end{remark}}
\def\bexer{\begin{exercise}}
\def\eexer{\end{exercise}}
\def\bfig{\begin{figure}}
\def\efig{\end{figure}}

\newtheorem{theorem}{Theorem}[section]
\newtheorem{prop}[theorem]{Proposition}
\newtheorem{lemma}[theorem]{Lemma}

\newtheorem{definition}[theorem]{Definition}

\newtheorem{example}[theorem]{Example}
\newtheorem{remark}[theorem]{Remark}

\begin{document}

\setlength{\baselineskip}{21pt}
\title{ Differentiating Matrix Functions \footnote{2010 \textit{MSC}: Primary 26B05, Secondary 15A60. }
   }
\author{Kelly Bickel 
\thanks{This research was partially supported by National Science Foundation grant DMS-0966845. }\\
Washington University\\
St. Louis, Missouri 63130 \\
kbickel@math.wustl.edu}
\date{}
\maketitle

\abstract{ Multivariate, real-valued functions on $\R^d$ induce matrix-valued functions on the space of $d$-tuples of $n \times n$ pairwise-commuting self-adjoint matrices. We examine the geometry of this space of matrices and conclude that the best notion of differentiation of these matrix functions is differentiation along curves. We prove that $C^1$ real-valued functions induces $C^1$ matrix functions and give a formula for the derivative. We also show that real-valued $C^m$ functions defined on open rectangles in $\R^2$ induce matrix functions that can be $m$-times continuously differentiated along $C^m$ curves. 

\bibliographystyle{plain}

\setlength{\parindent}{0in}
\section{Introduction}

Every real-valued function defined on $\R$ induces a matrix-valued function on the space of $n \times n$ self-adjoint matrices by acting on the spectrum of each matrix. Likewise, each  real-valued function $f$ defined on an open set $\Omega \subseteq \R^d$  induces a matrix-valued function $F$ on the space of $d$-tuples of $n \times n$ pairwise-commuting self-adjoint matrices with joint spectrum in $\Omega$. Let $S=(S^1, \dots, S^d)$ be such a $d$-tuple diagonalized by a unitary matrix $U$ as follows:
\beq S^r =
U \ \left( \begin{array}{ccc}
x^r_1 && \\
& \ddots & \\
&& x^r_n
\end{array} 
\right) U^* \qquad \forall \ 1 \le r \le d. 
\eeq
Denote the joint spectrum of $S$ by $\sigma(S):= \left \{ x_i=(x^1_i, \dots, x^d_i) : 1 \le i \le n \right \}$ and define 
\begin{eqnarray}
F(S):= U 
\left( \begin{array}{ccc}
f( x_1) && \\
&\ddots & \\
&& f( x_n)
\end{array} 
\right) \ U.^* \label{eqn1.1}
\end{eqnarray}

This paper will show that certain differentiability properties of the original function pass to the matrix function. Even for a one-variable function, this is nontrivial. Let $f \in C^1( \R, \R)$ and consider the simple case of differentiating the associated matrix function $F$ along a $C^1$ curve $S(t)$ of $n \times n$ self-adjoint matrices. At first glance, it seems reasonable to write $S(t)=U(t)D(t)U^*(t)$, for $U(t)$ unitary and $D(t)$ diagonal. Then $F(S(t)) = U(t)F(D(t))U^*(t)$ and we can differentiate using the product rule. \\

However, there is no guarantee that we can decompose $S(t)$ into its eigenvector and eigenvalue matrices so that the eigenvectors are even continuous. As demonstrated by the following example from \cite{rel37},  eigenvector behavior at points where distinct eigenvalues coalesce can be unpredictable. Specifically, let
\begin{eqnarray} \nonumber
S(t)\= e^{-\frac{1}{t^2}}\left(
\begin{array}{cc}
\cos( \frac{2}{t}) & \sin( \frac{2}{t}) \\
\\
\sin( \frac{2}{t}) & -\cos( \frac{2}{t})
\end{array} 
\right) \ \ \text{for t} \ne 0, \quad \text{ and } \quad S(0) =0.
\end{eqnarray}
For $t \ne 0$, the eigenvalues of $S(t)$ are $ \pm \ e^{-\frac{1}{t^2}}$ and their associated eigenvectors are
\begin{eqnarray} \nonumber
\pm \left(
\begin{array}{c}
\cos( \frac{1}{t}) \\
\\
\sin( \frac{1}{t}) \end{array} \right)
\text{ and }
\pm \left(
\begin{array}{c}
\sin( \frac{1}{t}) \\
\\
- \cos( \frac{1}{t}) \end{array} \right).
\end{eqnarray}
Thus, even an infinitely differentiable curve can have singularities in its eigenvectors. \\

The differentiability of matrix functions defined from one-variable functions is discussed frequently in the literature (see \cite{bh97}, \cite{don74}, \cite{horjoh91}). The most comprehensive result is by Brown and Vasudeva in \cite{brv00} , who prove that $m$-times continuously differentiable real functions induce $m$-times continuously Fr\'{e}chet differentiable matrix functions. \\

If a matrix function is defined using a real-valued function on $\R^d$ as in (\ref{eqn1.1}), its domain is the space of $d$-tuples of pairwise-commuting $n \times n$ self-adjoint matrices, denoted $CS^d_n$. For $d >1$, the space of $d$-tuples of $n \times n$ self-adjoint matrices is denoted $S^d_n$ and for $d=1$, is denoted $S_n.$\\

In Section 2, we analyze the geometry of $CS^d_n$ and conclude that the best notion of differentiability for functions on this space is differentiation along curves. If we fix $S$ in $CS^d_n$, Theorem \ref{thm2.1} characterizes the directions $\Delta$ in $S^d_n$ such that there is a $C^1$ curve $S(t)$ in $CS^d_n$ with $S(0)=S$ and $S'(0)= \Delta.$  In Theorem \ref{thm2.2}, we show that the joint eigenvalues of Lipschitz curves in $CS^d_n$ can be represented by Lipschitz functions. \\

In Section 3, we examine the differentiability properties of induced matrix functions. Specifically,  in Theorem \ref{thm3.1}, we  show that a $C^1$ function induces a matrix function that can be continuously differentiated along $C^1$ curves. We then calculate a formula for the derivative along curves and in Theorem \ref{thm3.2},  prove that it is continuous.\\

In Section 4, we consider higher-order differentiation.  With additional domain restrictions, in Theorem \ref{thm4.1}, we show that an induced matrix function is $m$-times continuously differentiable along $C^m$ curves. We also calculate a formula for the derivatives and in Theorem \ref{thm4.2}, show they are continuous.  In Section 5, we discuss several applications of the differentiability results.\\

There is an alternate approach for inducing a matrix function from a multivariate function; the $d$ matrices $S^1, \dots, S^d$ are viewed as operators on Hilbert spaces $H^1, \dots, H^d$ and $F(S)$ is viewed as an operator on $H^1 \otimes \dots \otimes H^d$. Brown and Vasudeva generalize their one-variable result to these matrix functions in \cite{brv00}.\\

Before proceeding, I  would like to thank John McCarthy for his guidance during this research and the referees for their many useful suggestions.

\section{The Geometry of $CS^d_n$}
Let $S=(S^1, \dots,  S^d)$ be in  $CS^d_n$  (or  $S^d_n$) and let  $x_i =(x^1_i, \dots, x^d_i)$ be in $\sigma(S).$ Define
$$ \| S \| := \max_{1 \le r \le d} \| S^r\| \ \ \ \text{and} \ \ \   \| x_i\| := \max_{1 \le r \le d} |x^r_i|,$$
where $\|S^r\|$ is the usual operator norm.  Observe that $CS^d_n$ is not a linear space; if $A$ and $B$ are pairwise-commuting $d$-tuples, the sum $A+B$ need not pairwise commute. Thus, neither the Fr\'echet nor G\^{a}teaux derivatives can be defined for functions on $CS^d_n$ because both require the function to be defined on linear sets around each point. \\

Recall that $CS^d_n$ is the set of elements $S \in S^d_n$ with $[ S^r, S^s]=0$ for all $ 1 \le r,s \le d.$  Thus, $CS^d_n$ is the zero set of the polynomials associated with $\frac{d(d-1)}{2}$ commutator operations and so is an algebraic variety. A result by Whitney \cite{whi65} says every algebraic variety can be decomposed into submanifolds that fit together `regularly' and whose tangent spaces fit together `regularly.' For a manifold $N$, let $TN$ denote the tangent space of $N$ and let $T_SN$ denote the tangent space based at a point $S$ in $N$. To make Whitney's conditions more precise, we need the following definition:
\begin{definition}
A \textbf{stratification} of $X$ is a locally finite partition $Z$ of $X$ such that
\beq
&(i) & \text{ Each piece $ M_{\alpha} \in Z $ is a smooth submanifold of $X.$} \\
&(ii) & \text{The frontier of each piece $ \overline{M}_{\alpha} \ \backslash M_{\alpha} $ is either trivial or a union of other pieces.}
\eeq
Then $X$ is called a \textbf{stratified space with stratification $Z$}.
\end{definition}
\begin{example} \label{ex2.1} Consider $CS^2_2$, the space of pairs of self-adjoint, commuting $2 \times 2$ matrices. In the following definitions, $a, b, c, d \in \R$. Define
\beq M_1 &: =&  \Big \{  \left(  U \left( 
\begin{array}{cc}
a & 0 \\ 
0 & b
\end{array} 
\right)U^* , \ U\left( 
\begin{array}{cc}
c & 0 \\ 
0 & d
\end{array} 
\right)U^*  
\right)
: \ U \in S_2 \text{ is unitary, }  a\ne b, \  c \ne d 
\Big \} \eeq
\beq M_2 &: = & \Big \{  \left(  \left( 
\begin{array}{cc}
a & 0 \\ 
0 & a
\end{array} 
\right) , \ U\left( 
\begin{array}{cc}
c & 0 \\ 
0 & d
\end{array} 
\right)U^*  
\right)
: \ U \in S_2 \text{ is unitary, }  \  c \ne d 
\Big \} \hspace{.5in} \\
 M_3 &: =&  \Big \{  \left(  U \left( 
\begin{array}{cc}
a & 0 \\ 
0 & b
\end{array} 
\right)U^* , \ \left( 
\begin{array}{cc}
c & 0 \\ 
0 & c
\end{array} 
\right)  
\right)
: \ U \in S_2 \text{ is unitary, } \  a \ne b 
\Big \}\\
M_4 &: =&  \Big \{  \left(   \left( 
\begin{array}{cc}
a & 0 \\ 
0 & a
\end{array} 
\right) , \left( 
\begin{array}{cc}
c & 0 \\ 
0 & c
\end{array} 
\right)  
\right) 
\Big \} 
\end{eqnarray*}
It is clear that $CS^2_2 = \cup M_i$. Moreover, each $M_i$ is a manifold and $\overline{M}_i \ \backslash M_i$ is either trivial or a union of other $M_j.$ Thus, the partition $\{M_i\}$ is a stratification of $CS^2_2.$
\end{example}
In general, a decomposition of $CS^d_n$ into pieces will be related to the number and multiplicity of the repeated joint eigenvalues of the elements of $CS^d_n$. \\

Whitney's result  says $CS^d_n$ has a stratification $Z$ with further regularity. Specifically, let $\{M_{\alpha} \}$ denote the pieces of $Z$ and define $TCS^d_n:= \cup TM_{\alpha}$. Then, $TCS^d_n$ is also a stratified space, and we call $Z$ a Whitney stratification of $CS^d_n$. Given a function $F: CS^d_n \rightarrow S_n$, one type of derivative is a  map
\beq  DF: TCS^d_n \rightarrow TS_n \ \ \text{ such that} \ \ \  DF|_{TM_{\alpha}}: TM_{\alpha} \rightarrow TS_n \eeq 
is the usual differential map for each $M_{\alpha}.$ In Theorem \ref{thm3.3}, we analyze such maps. However, these differential maps cannot be easily generalized to analyze higher-order differentiation. Furthermore, the space $TCS^d_n$ will only contain a subset of the vectors tangent to $CS^d_n$. Example \ref{ex2.2}  will show that strict containment often occurs.\\

To retain information about all tangent vectors, we will mostly study differentiation along differentiable curves. We first determine which $\De$ in $S^d_n$ are vectors tangent to $CS^d_n$ at a given point $S.$ This is equivalent to the following question:
$$ \text{ Is there a } C^1 \text{ curve } S(t) \text{ in } CS^d_n \text{ with } S(0)=S \text{ and } S'(0) = \De?$$
For an element $S \in CS^d_n$ with distinct joint eigenvalues, Agler, McCarthy, and Young in \cite{mcc10} gave necessary and sufficient conditions on $S$ and $\De$ for such a $C^1$ curve to exist. We extend their result to an arbitrary element $S$. Fix $S \in CS^d_n$ and $\De \in S^d_n$. Let $U$ be a unitary matrix diagonalizing each component of $S$ such that the repeated joint eigenvalues appear consecutively. Renumbering the $x_i$'s if necessary, define
\begin{eqnarray}
\label{eq2.1}
D^r &:=& U^* S^rU \= \left(
\begin{array}{ccc}
x_1^r && \\
&\ddots & \\
&& x_n^r 
\end{array} 
\right) 
\qquad \forall \ 1 \leq r \leq d.
\end{eqnarray}
 For each $r$, define the two matrices
\begin{eqnarray} \Gamma^r &:=& U^* \ \De^r \ U \nonumber \\ 
&& \nonumber \\
\label{eq2}
\tilde{\Gamma}^r_{ij} 
&:=& \left\{
\begin{array}{ll}
\Gamma^r_{ij} & \mbox{if\ } x_i =x_j \\ 
0 & \mbox{otherwise.\ } 
\end{array} \right.
\end{eqnarray}
Then $\tilde{\Gamma}^r$ is a block diagonal matrix. Each block corresponds to a distinct joint eigenvalue of $S$ and has dimension equal to the multiplicity of that eigenvalue.

\bt \label{thm2.1}
Let $S \in CS^d_n$ and $\De \in S^d_n.$ There exists a $C^1$ curve $S(t)$ in $CS^d_n$ with $S(0) =S$ and $S'(0) = \De$ iff $$ \ \left[ D^r, \Gamma^s \right] = \left[D^s, \Gamma^r \right] \text{ and } \big [ \tilde{ \Gamma}^r, \tilde{ \Gamma}^s \big ] =0 \qquad \forall \ 1 \le r,s \le d.$$ \et
\bp $( \Rightarrow )$ Assume $S(t) $ is a $C^1$ curve in $CS^d_n$ with $S(0)=S$ and $S'(0) = \De.$ Define 
$$R(t):= U^*\ S(t) \ U,$$ 
where $U$ diagonalizes $S$ as in (\ref{eq2.1}). Then $R(t)$ is a $C^1$ curve in $CS^d_n$ with $S(0)=D$ and $S'(0) = \Gamma$. We will first prove that
$$\left[ D^r, \Gamma^s \right] = \left[ D^s, \Gamma^r \right] \text{ and } \left[ \Gamma^r, \Gamma^s \right]_{ij} =0 \quad \forall \ 1 \le r,s \le d \text{ and } (ij) \text{ such that } x_i=x_j.$$
We will use those commutativity results to conclude
$$\big [ \tilde{ \Gamma}^r, \tilde{ \Gamma}^s \big ] =0 \qquad \forall \ 1 \le r,s \le d.$$
Since $R(t)$ is $ C^1$ in a neighborhood of $t=0$, we can write
$$R^r(t) = D^r + \Gamma^rt+h^r(t)\qquad \forall \ 1 \le r \le d,$$
where $| h^r(t)_{ij} |=o( |t| )$ for $ 1 \le i,j \le n.$ For each pair $r$ and $s$, the pairwise-commutativity of $R(t)$ implies
\begin{eqnarray} 0 \ &=& [R^r(t), R^s(t) ] \nonumber \\
&=& [ \ D^r + \Gamma^rt+h^r(t), \ D^s + \Gamma^st+h^s(t) \ ] \nonumber \\
&=& \big( [D^r,h^s(t) \ ]+ [ h^r(t), D^s]+ [ h^r(t), \ h^s(t)] \big) \nonumber \\
&& \ \ + \ \big ( \ [D^r, \Gamma^s] + [\Gamma^r,D^s] + [\Gamma^r, h^s(t)] + [h^r(t), \Gamma^s] \ \big ) t \nonumber \\
& & \ \ + \ \ \ [ \Gamma^r, \Gamma^s] \ t^2 ,\label{eq2.2}\end{eqnarray}
where the term $[D^r,D^s]$ was omitted because it vanishes. Fix $t \ne 0$ and  divide each term in (\ref{eq2.2}) by $t$. Letting $t$ tend towards zero yields
\begin{eqnarray} 0 &=& [D^r, \Gamma^s] - [D^s, \Gamma^r ]. \label{eq2.2.1}\end{eqnarray}
Choose $i$ and $j$ such that $x_i = x_j$. Then, the $ij^{th}$ entry of (\ref{eq2.2}) reduces to
$$ 0 = [ h^r(t),h^s(t)]_{ij} + \left( \ [\Gamma^r, h^s(t)]_{ij} - [\Gamma^s, h^r(t)]_{ij} \ \right )t + [\Gamma^r, \Gamma^s]_{ij}t^2 .$$
Fix $t \ne 0$ and divide both sides by $t^2.$ Letting $t$ tend towards zero yields
\begin{eqnarray} \label{eq2.3.2} 0 = [\Gamma^r, \Gamma^s]_{ij}. \end{eqnarray}
Fix $ r$ and $s$ with $ 1 \le r,s\le d.$ Since $ \tilde{\Gamma}^r$ and $ \tilde{\Gamma}^s$ are block diagonal matrices with blocks corresponding to the distinct joint eigenvalues of $S,$ it follows that $\tilde{\Gamma}^r\tilde{\Gamma}^s$ and $\tilde{\Gamma}^s\tilde{\Gamma}^r$ are also such block diagonal matrices. Thus, if $i$ and $j$ are such that $x_i \ne x_j,$
\beq \big [ \ \tilde{\Gamma}^r, \ \tilde{\Gamma}^s \big ]_{ij} &=&\big ( \ \tilde{\Gamma}^r\tilde{\Gamma}^s - \tilde{\Gamma}^s\tilde{\Gamma}^r \ \big )_{ij} = 0.\eeq
Now, fix $i$ and $j$ such that $x_i=x_j$. 
By the definition of $\tilde{\Gamma}$,
\beq \big [\tilde{\Gamma}^r, \tilde{\Gamma}^s \big ]_{ij} &=& \ \ \ \sum_{k=1}^n \ \ \ \tilde{\Gamma}^r_{ik}\tilde{ \Gamma}^s_{kj} - \tilde{ \Gamma}^s_{ik}\tilde{\Gamma}^r_{kj} \\
&=& \sum_{ \left \{ k : x_k = x_i \right \} } \Gamma^r_{ik}\Gamma^s_{kj} - \Gamma^s_{ik} \Gamma^r_{kj} \\ 
&=& - \sum_{ \left \{ k : x_k \ne x_i \right \} } \Gamma^r_{ik}\Gamma^s_{kj} - \Gamma^s_{ik} \Gamma^r_{kj}, \eeq
where the last equality uses (\ref{eq2.3.2}). Thus, it suffices to show that if $x_k \ne x_i$, 
$$ \Gamma^r_{ik} \Gamma^s_{kj} - \Gamma^s_{ik} \Gamma^r_{kj} =0.$$
Assume $x_k \ne x_i,$ and fix $q$ with $x^q_k \ne x^q_i.$ Apply (\ref{eq2.2.1}) to pairs $r,q$ and $s,q$ to get
 $$[ D^q, \Gamma^r] = [D^r, \Gamma^q] \ \ \text{ and } \ \ [ D^q, \Gamma^s] = [D^s, \Gamma^q].$$
Restricting to the $ik^{th}$ and $kj^{th}$ entries of the previous two equations yields
\begin{eqnarray} \label{eq2.3.1}
\begin{array}{ccc}
\Gamma^r_{ik} (x^q_i-x^q_k) &=& \Gamma^q_{ik}(x^r_i-x^r_k) \ \ \ \ \ \Gamma^r_{kj} (x^q_k-x^q_j) \ \ = \ \ \Gamma^q_{kj}(x^r_k-x^r_j)  \\ 
&&\\
\Gamma^s_{ik} (x^q_i-x^q_k) &=& \Gamma^q_{ik}(x^s_i-x^s_k)  \ \ \ \ \ \Gamma^s_{kj} (x^q_k-x^q_j)\ \  = \ \ \Gamma^q_{kj}(x^s_k-x^s_j).  
\end{array} \end{eqnarray}
Since $x_i =x_j$ and $x^q_k \ne x^q_i$, we can replace all the $x_j$'s with $x_i$'s in $(\ref{eq2.3.1} )$ and solve for the $\Gamma^r$ and $\Gamma^s$ entries. Using these relations gives
$$ \Gamma^r_{ik} \Gamma^s_{kj} - \Gamma^s_{ik} \Gamma^r_{kj} = \dfrac {\Gamma^q_{ik}(x^r_i-x^r_k) \Gamma^q_{kj}(x^s_i-x^s_k)}{ (x^q_i-x^q_k)^2}-
\dfrac {\Gamma^q_{ik}(x^s_i-x^s_k) \Gamma^q_{kj}(x^r_i-x^r_k)}{(x^q_i-x^q_k)^2 }=0, $$
as desired. Thus, $ [ \tilde{\Gamma}^r, \tilde{\Gamma}^s ] =0.$ \\

$( \Leftarrow)$ Fix $S \in CS^d_n$ and $\De \in S^d_n$ and let $U$, $D$, and $\Gamma$ be as in the discussion preceding Theorem \ref{thm2.1}. Assume
\begin{eqnarray} \label{eq2.4} \left[ D^r, \Gamma^s \right] = \left[ D^s, \Gamma^r \right] \text{ and } \big [ \tilde{ \Gamma}^r, \tilde{ \Gamma}^s \big ] = 0 \qquad \forall \ 1 \le r,s \le d. \end{eqnarray}
Define a skew-Hermitian matrix $Y$ as follows:
\beq
Y_{ij} 
:= \left\{
\begin{array}{ll}
\tfrac{\Gamma^q_{ij}}{x^q_j -x^q_i} & \mbox{if\ } x_i \ne x_j \\ 
0 & \mbox{otherwise,\ } \\
\end{array} \right. \\
\eeq
where the $q$ is chosen so that $x^q_i -x^q_j \ne 0$. Observe that $Y$ is independent of $q$ because the $ij^{th}$ entry of the first equation in (\ref{eq2.4}) is
$$ \Gamma^s_{ij} (x^r_i-x^r_j) = \Gamma^r_{ij} (x^s_i-x^s_j).$$
Now, define the curve $S(t)$ by
$$S^r(t) := U e^{Yt} \big [ D^r + t \tilde{\Gamma}^r \big ] e^{-Yt} U^* \qquad \forall \ 1 \le r \le d.$$
Then, $S(t)$ is continuously differentiable. Because $Y$ is skew-Hermitian, $e^{Yt}$ is unitary. Since $D^r$ and $\tilde{\Gamma}^r$ are self-adjoint, $S(t) \in S^d_n.$ By a simple calculation using (\ref{eq2.4}), 
$$\left[ S^r(t), S^s(t) \right] =0 \quad \forall \ 1 \le r,s \le d.$$
Thus, $S(t) \in CS^d_n.$ By definition, $S(0) = S.$ For each $r$, 
$$(S^r)'(t) =U \left( Ye^{Yt} \big [ D^r + t \tilde{\Gamma}^r \big ] e^{-Yt}+ e^{Yt} \big [ \tilde{\Gamma}^r \big ] e^{-Yt} - e^{Yt} \big [ D^r + t \tilde{\Gamma}^r \big ]Y e^{-Yt} \right) U^*,$$
so that
$$ (S^r)'(0) = U \left( \left[ Y, D^r \right] +\tilde{\Gamma}^r \right) U^*= \De^r.$$ 
Thus, $S'(0) = \De,$ and $S(t)$ is the desired curve.\ep
\vspace{.1in}
\begin{example} \label{ex2.2} Let $I \in CS^d_n$ be the identity element. By Theorem \ref{thm2.1}, there is a continuously differentiable curve $S(t)$ in $CS^d_n$ with
$$ S(0)=I \text{ and } S'(0) = \De \ \ \text{ if and only if }\ \ \De \in CS^d_n.$$ 
Thus, the set of vectors tangent to $CS^d_n$ at $I$ is $CS^d_n.$ For a Whitney stratification of $CS^d_n$ and piece $M_{\alpha}$ containing $I$, the tangent space $T_IM_{\alpha}$  is linear. Since $CS^d_n$ is not linear, $T_IM_{\alpha}$ is a strict subset of the set of tangent vectors at $I$. 
\end{example}  

The conditions of Theorem \ref{thm2.1} actually imply that if $S$ in $CS^d_n$ has any repeated joint eigenvalues, the set of vectors tangent to $CS^d_n$ at $S$ is not a linear set.  Then, for any Whitney stratification of $CS^d_n$ and piece $M_{\alpha}$ containing $S$, the tangent space $T_SM_{\alpha}$ is a strict subset of the vectors tangent to $CS^d_n$ at $S.$ We will thus focus on differentiation along curves rather than differential maps. \\

To evaluate an induced matrix function along a curve in $CS^d_n$, we apply the original function to curve's joint eigenvalues. We are therefore interested in the behavior of the joint eigenvalues of curves in $CS^d_n$. \\

If $S(t)$ is a continuous curve in $S_n,$ a result by Rellich in \cite{rel37} and \cite{rel37b} states that the eigenvalues of $S(t)$ can be represented by $n$ continuous functions. A succinct proof is given by Kato in \cite[pg 107-10]{kat66}. With slight modification, the arguments  show that the eigenvalues of a Lipschitz curve in $S_n$ can be represented by Lipschitz functions. These results generalize as follows:

\bt \label{thm2.2} Given a Lipschitz curve $S(t)$ in $CS^d_n$ defined on an interval $I$, there exist Lipschitz functions $x_1(t),..., x_n(t) : I \rightarrow \R^d$ with $\sigma(S(t)) = \left \{ x_i(t) : 1 \le i \le n \right \}.$ \et

\bp As the proof is a technical but straightforward  modification of the one-variable case, it is left as an exercise. \ep

 Theorem \ref{thm2.2} provides a specific ordering of the eigenvalues of $S(t)$ at each $t.$ This ordering may differ from the one in (\ref{eq2.1}), where joint eigenvalues appear consecutively.  However, Theorem \ref{thm2.2} implies that the eigenvalues of a Lipschitz curve $S(t)$ are Lipschitz as an unordered $n$-tuple. Specifically, fix $t^*$ and denote the eigenvalues of $S(t^*)$ by $\{x_i : 1 \le i \le n \}$.  Then, for $t$ near $t^*$, there is a constant $c$ such that
$$ \min \bigg ( \max_{1 \le i \le n} \| x_i - x_i(t) \| \bigg)  \le c|t^* -t |,$$ 
where the minimum is taking over all reorderings of the $\{x_i\}.$  If we require that eigenvalues are ordered as in (\ref{eq2.1}), we will use Theorem \ref{thm2.2} to conclude that the eigenvalues are Lipschitz as an unordered $n$-tuple. \\

\section{Differentiating Matrix Functions}
Recall that every real-valued function defined on an open set $\Omega \subseteq R^d$ induces a matrix function as in (\ref{eqn1.1}). We denote its domain, the space of $d$-tuples of pairwise-commuting $n \times n$ self-adjoint matrices with spectrum in $\Omega$, by $CS^d_n(\Omega)$. \\

If the original function is continuous, the matrix function is as well. Specifically, Horn and Johnson proved in \cite[pg 387-9]{horjoh91} that a one-variable polynomial induces a continuous matrix polynomial. The arguments generalize easily to multivariate polynomials, and approximation arguments imply that the matrix function of a continuous function is continuous. We now consider differentiability and prove:
 
\bt \label{thm3.1} Let $S(t)$ be a $C^1$ curve in $CS^d_n$ defined on an interval $I$ and let $\Omega$ be an open set in $\R^d$ with $ \sigma(S(t)) \subset \Omega.$ If $f \in C^1(\Omega, \R)$, then
\beq 
&(i)& \tfrac{d}{dt} F(S(t)) |_{t=t^*} \text{ exists for all } t^* \in I. \\
&(ii)& \text{If } T(t) \text{ is another } C^1 \text{ curve in } CS^d_n \text{ with } T(0)=S(t^*) \text{ and } T'(0) = S'(t^*), \text{ then} \eeq
$$ \tfrac{d}{dt} F(T(t)) |_{t=0} = \tfrac{d}{dt} F(S(t)) |_{t=t^*}. $$ \et

\vspace{.1in}

We say an open set $\Omega \subset \R^d$ is a \textbf{rectangle} if $\Omega =I^1 \times \dots \times I^d,$ and an open set $\tilde{\Omega} \subset \C^d$ is a \textbf{complex rectangle} if $\tilde{\Omega} = (I^1+iJ^1) \times \dots \times (I^d+iJ^d),$ where each $I^r$ and $J^r$ is an open interval in $\R$. Before proving Theorem \ref{thm3.1}, we assume $f$ is real-analytic and prove Proposition \ref{prop3.1}.  See \cite{horjoh91} for the one-variable case.  

\bprop \label{prop3.1} Let $S(t)$ be a $C^1$ curve in $CS^d_n$ defined on an interval $I$. Let $\Omega$ be an open rectangle in $\R^d$ with $ \sigma(S(t)) \subset \Omega.$ If $f$ is a real-analytic function on $\Omega$, then
$$ \tfrac{d}{dt} F(S(t)) |_{t=t^*} \text{ exists and is continuous as a function of } t^* \text{ on } I.$$ \eprop 

\vspace{.1in}

The proof of Proposition \ref{prop3.1} requires the following two lemmas.

\bl  \label{re3.1} Let $\Omega$ be an open rectangle in $\R^d$ and let $S \in CS^d_n$ with $ \sigma(S) \subset \Omega.$ Each real-analytic function  on $\Omega$ can be extended to an analytic function defined on a complex rectangle $\tilde{\Omega}$ such that $\sigma(S)$ is in $\tilde{\Omega}.$ \el
\bp The result follows from basic properties of complex functions. It should be noted that $\tilde{\Omega}$ need not contain $\Omega$.\ep

\bl \label{lem3.1} Let $\tilde{\Omega}$ be an open rectangle in $ \C^d$ and let $S \in CS^d_n$ with $\sigma(S) \subset \tilde{\Omega}.$ If $f$ is an analytic function on $\tilde{\Omega},$ then
\begin{eqnarray} \nonumber \dis F(S) = \frac{1}{(2 \pi i)^d} \int_{C^d} \dots \int_{C^1} f( \zeta^1, \dots, \zeta^d) (\zeta^1I -S^1)^{-1} \dots (\zeta^dI -S^d)^{-1} \ d\zeta^1\dots d\zeta^d, \end{eqnarray}
where $C^r$ is a rectifiable curve strictly containing $\sigma(S^r)$, and $C^1 \times \dots \times C^d \subset \tilde{\Omega}$. \el
\bp Horn and Johnson prove the formula for a one-variable function in \cite[ pg 427]{horjoh91}. Their derivation generalizes easily to multivariate functions.\ep

\vspace{.1in}

Proof of Proposition \ref{prop3.1}: \\  
For ease of notation, assume $d=2$ and define
$$ R^r(t) :=( \zeta^r I -S^r(t))^{-1} \qquad \forall \ 1 \le r \le 2,$$
where $\zeta^r$ is in the resolvent of $S^r(t).$ Fix $t_0 \in I$ and extend $f$ to an analytic function on a complex open rectangle $\tilde{\Omega}$ containing $\sigma(S(t_0)).$ Choose rectifiable curves $C^1$ and $C^2$ such that $C^1 \times C^2 \subset \tilde{\Omega}$ and each $C^r$ strictly encloses the eigenvalues of $S^r(t_0)$.  By Theorem \ref{thm2.2}, the joint eigenvalues of $S(t)$ are continuous and by Lemma \ref{lem3.1},
\beq \dis F(S(t)) = \frac{1}{(2 \pi i)^2} \int_{C^2} \int_{C^1} f( \zeta^1, \zeta^2)\ R^1(t) \ R^2(t) \ d\zeta^1 d\zeta^2, \eeq 
for $t$ sufficiently close to $t_0$. Direct calculation gives
$$ \tfrac{d}{dt} R^r(t) |_{t=t^*}= R^r(t^*) \ (S^r)'(t^*) \ R^r(t^*) \quad \text{ for } 1 \le r \le 2 \text{ and } t^* \text{ near } t_0.$$
It can be easily shown that, for $t^*$ sufficiently close to $t_0$, we can interchange integration and differentiation to yield
\begin{eqnarray} 
\dis \tfrac{d}{dt} F(S(t)) |_{t=t^*} &=&  \frac{1}{(2 \pi i)^2} \int_{C_2} \int_{C_1} f( \zeta^1, \zeta^2) \ \frac{d}{dt} \bigg ( R^1(t)R^2(t) \bigg ) \Big |_{t=t^*} \ d\zeta^1 d\zeta^2 \nonumber \\
&& \nonumber \\
& =&\frac{1}{(2 \pi i)^2} \int_{C_2} \int_{C_1} f( \zeta^1, \zeta^2) \bigg ( R^1(t^*) \ (S^1)'(t^*) \ R^1(t^*)R^2(t^*) \nonumber \\
&& \nonumber \\
&& \ \ \ \ \ + \ R^1(t^*)R^2(t^*) \ (S^2)'(t^*) \ R^2(t^*) \bigg ) d\zeta^1 d\zeta^2. \label{eq3.1}
\end{eqnarray}
As each $(S^r)'(t)$ is continuous and all other terms in (\ref{eq3.1}) are uniformly bounded near $t_0$, we get $\frac{d}{dt} F(S(t)) |_{t=t^*}$ is continuous at $t^*=t_0$. \ep 

\vspace{.1in}

Proof of Theorem \ref{thm3.1}: \\
Observe that the theorem holds for polynomials:  $(i)$ follows from Proposition \ref{prop3.1} and $(ii)$ follows from the formula in (\ref{eq3.1}). Fix $t^* \in I$. Let $f$ be an arbitrary $C^1$ function and let $p$ be a polynomial that agrees with $f$ to first order on $\sigma(S(t^*)).$ \\

By Theorem \ref{thm2.2}, there are Lipschitz maps $x_i(t):= (x_i^1(t), \dots, x_i^d(t)),$ for $1 \le i \le n,$  representing $\sigma(S(t))$ on $I$. From the multivariate Mean Value Theorem, we have
\begin{eqnarray} 
\|(F-P)(S(t)) \| &\= & \dis \max_{i} | (f-p)(x_i(t)) | \nonumber \\
&=& \dis \max_{i} | (f-p)(x_i(t))-(f-p)(x_i(t^*)) \big | \nonumber \\
&=& \dis \max_{i} \big | \nabla (f-p)(x^*_i(t)) \cdot \left( x_i(t) -x_i(t^*) \right) \big | \nonumber \\
& \le & \dis \max_{i} \dis \sum_{r=1}^d \big | \left( \tfrac{ \partial f}{\partial x^r} - \tfrac{ \partial p}{\partial x^r} \right) ( x^*_i(t)) \big |  | \ x^r_i(t) -x^r_i(t^*)  |, \label{eq3.2.1}
\end{eqnarray}
where $x^*_i(t)$ is on the line connecting $x_i(t)$ and $x_i(t^*)$ in $\R^d$. For $t$ near $t^*$, continuity implies $x^*_i(t) \in \Omega.$ As $f$ and $p$ agree to first order on $\sigma(S(t^*)),$ from (\ref{eq3.2.1}), we have
$$\|(F-P)(S(t))\| = o(|t- t^*|).$$ Hence
$$ \Big \| \frac{ F(S(t))-F(S(t^*))}{t-t^*} - \frac{P(S(t)) - P(S(t^*))}{t-t^*} \Big \| \rightarrow 0 \qquad \text{ as } t \rightarrow t^*. $$ 
Therefore, 
\beq \tfrac{d}{dt} F(S(t)) |_{t=t^*} \ \ \ \text{ exists and equals } \ \ \ \tfrac{d}{dt} P(S(t)) |_{t=t^*}. \eeq 
Applying the same argument to $F(T(t))$ at $t=0$ gives
\beq \tfrac{d}{dt} F(T(t)) |_{t=0} \ \ \ \text{ exists and equals } \ \ \ \tfrac{d}{dt} P(T(t)) |_{t=0}. \eeq 
As $(ii)$ holds for $P(t)$,  we must have 
\beq\hspace{1.75in} \tfrac{d}{dt} F(T(t)) |_{t=0} = \tfrac{d}{dt} F(S(t)) |_{t=t^*}. \hspace{1.65in} \text{\hfill {$\Box$}}\eeq 

\vspace{.1in}

In the following proposition, we calculate an explicit formula for the derivative.

\vspace{.1in}
\bprop \label{prop3.2} Let $S(t)$ be a $C^1$ curve in $CS^d_n$ defined on an interval $I$ and let $t^* \in I$. Let $\Omega$ be an open set in $\R^d$ with $ \sigma(S(t)) \subset \Omega$ and let $f \in C^1(\Omega,\R)$. Then,
$$ \tfrac{d}{dt} F(S(t)) |_{t=t*} = U \bigg ( \sum_{r=1}^d \ \tilde{\Gamma}^r \tfrac{ \partial f}{\partial x^r}(D) + \left[ Y, F(D) \right] \bigg ) U^*,$$
where $U$ diagonalizes $S(t^*)$ as in $(\ref{eq2.1})$ and the other matrices are as follows:
\begin{eqnarray*}
\nonumber D^r &:=& U^* \ \big [ S^r(t^*) \big ] \ U  \ \ \ \ \ \ \ \ \ \Gamma^r := \ U^* \ \big [ (S^{r})' (t^*) \big ] \ U \ \\ 
\nonumber \tilde{\Gamma}^r_{ij} 
&:=& \left\{
\begin{array}{ll}
\Gamma^r_{ij} & \mbox{if\ } x_i =x_j \\ 
0 & \mbox{otherwise} 
\end{array} \right. 
\ \ \  Y_{ij} 
:= \left\{
\begin{array}{ll}
\tfrac{\Gamma^q_{ij}}{x^q_j -x^q_i} & \mbox{if\ } x_i \ne x_j \\ 
0 & \mbox{otherwise, \ } \end{array} \right. \end{eqnarray*} 
where the joint eigenvalues of $S(t^*)$ are given by $\left \{ x_i = (x^1_i, \dots ,x^d_i) : \ 1 \le i \le n \right \}$ and $q$ is chosen so  $x^q_j -x^q_i \ne 0.$  
\eprop

\bp  Let $t^* \in I$ and define the $C^1$ curve $T(t)$ by
$$T^r(t) := U \ e^{Yt} \big [ D^r + t \tilde{\Gamma}^r \big ] e^{-Yt} \ U^* \qquad \forall \ 1 \le r \le d.$$
Then, $T(t)$ is the curve defined in the proof of Theorem \ref{thm2.1} for $S:=S(t^*)$ and $\De:=S'(t^*)$. It is immediate that $T(t) \in CS^d_n$, $T(0)=S(t^*)$, and $T'(0) = S'(t^*)$. By Theorem \ref{thm3.1}, it now suffices to calculate $\frac{d}{dt} F(T(t)) |_{t=0}.$  First, we diagonalize each $ D^r + t \tilde{\Gamma}^r .$ Let $p$ be the number of distinct joint eigenvalues of $S(t^*)$. By definition,
\begin{eqnarray} \nonumber \tilde{\Gamma}^r\= \left(
\begin{array}{ccc}
\Gamma^r_1 && \\
&\ddots& \\
&& \Gamma^r_p 
\end{array} 
\right) \qquad \forall \ 1 \le r \le d,
\end{eqnarray}
where each $\Gamma^r_l$ is a $k_l \times k_l$ self-adjoint matrix corresponding to a distinct joint eigenvalue of $S$ with multiplicity $k_l.$ It follows from Theorem \ref{thm2.1} that
$$ \big [ \tilde{ \Gamma}^r, \tilde{ \Gamma}^s \big ] =0, \ \text{ which implies: } \ \big [\Gamma^r_l, \Gamma^s_l \big ] =0 \qquad \forall \ 1 \le r,s \le d \text{ and } 1 \le l \le p.$$
Thus, for each $l$, there is a $k_l \times k_l$ unitary matrix $V_l$ such that $V_l$ diagonalizes each $\Gamma^r_l$. Let $V$ be the $n \times n$ block diagonal matrix with blocks given by $V_1, \dots, V_p.$ Then, $V$ is a unitary matrix that diagonalizes each $\tilde{\Gamma}^r$. By the diagonalization in (\ref{eq2.1}), the joint eigenvalues of $D$ are positioned so that
\begin{eqnarray} \label{eq3.3}
D^r\= \left(
\begin{array}{ccc}
c^r_1 I_{k_1} && \\
&\ddots& \\
&& c^r_p I_{k_p}
\end{array} 
\right) \qquad \forall \ 1 \le r \le d,
\end{eqnarray}
where $I_{k_l}$ is the $k_l \times k_l$ identity matrix and $c^r_l$ is a constant. Equation (\ref{eq3.3}) shows that conjugation by $V$ will not affect $D^r$. Define the diagonal matrix
$$\Lambda^r:= V^* \ \tilde{\Gamma}^r \ V \qquad \forall \ 1 \le r \le d,$$ 
and rewrite $T(t)$ as follows
\beq T^r(t) &=& U \ e^{Yt} V \left[ D^r + t \Lambda^r \right]V^{*} e^{-Yt} \ U^* \qquad \forall \ 1 \leq r \leq d. \eeq
Now we directly calculate $F(T(t))$ and $\frac{d}{dt} F(T(t)) |_{t=0}$
\begin{eqnarray*} F(T(t)) &=& U e^{Yt} V \ F \left( D^1 + t \Lambda^1 , \dots, D^d + t \Lambda^d \right) \ V^{*} e^{-Yt} U^*  \\
&=&  U e^{Yt}  V \ \bigg ( F(D) + t \dis \sum_{r=1}^d \Lambda^r \tfrac{ \partial f}{\partial x^r}(D) +o(|t|  ) \bigg ) \ V^{*}  e^{-Yt}  U^* , \label{eq3.4} \end{eqnarray*} 
where$\frac{ \partial f}{\partial x^r}(D) $ is defined by
\beq \tfrac{ \partial f}{\partial x^r}(D) := \left ( \begin{array}{ccc}
\frac{ \partial f}{\partial x^r}(x_1) && \\
& \ddots & \\
&& \frac{ \partial f}{\partial x^r}(x_n) \end{array} \right ) \qquad \forall \ 1 \le r \le d, \eeq 
and the first-order approximation of $F$ follows from the approximation of $f$ on each diagonal entry of the $d$-tuple of diagonal matrices. Differentiating $F(T(t))$ and setting $t=0$ gives
\beq  \tfrac{d}{dt} F(T(t))  |_{t=0} &=& U \Big( \sum_{r=1}^d V \ \Lambda^r \tfrac{ \partial f}{\partial x^r}(D) V^{*} + \left[ Y, V F(D) V^{*} \right] \Big) U^* \nonumber \\
&=& U \Big( \sum_{r=1}^d \ \tilde{\Gamma}^r \tfrac{ \partial f}{\partial x^r}(D) + \left[ Y, F(D) \right] \Big) U^*, \eeq
where conjugation by $V$ leaves $F(D)$ and each $\frac{ \partial f}{\partial x^r}(D)$ unchanged because those matrices have decompositions akin to that of $D^r$ in (\ref{eq3.3}). \ep

\vspace{.1in}

We now prove that the derivative calculated in Proposition \ref{prop3.2} is continuous in $t^*$. 

\bt \label{thm3.2}Let $S(t)$ be a $C^1$ curve in $CS^d_n$ defined on an interval $I$. Let $\Omega$ be an open set in $\R^d$ with $ \sigma(S(t)) \subset \Omega.$ If $f \in C^1(\Omega, \R)$, then 
$$ \tfrac{d}{dt} F(S(t)) |_{t=t^*} \text{ is continuous as a function of } t^* \text{on } I.$$ \et

For the proof, we will require the following lemma:

\bl \label{lem3.4} Let $S(t)$ be a $C^1$ curve in $CS^d_n$ defined on an interval $I$. Let $\Omega$ be an open, convex set in $\R^d$ with $ \sigma(S(t)) \subset \Omega.$ If $f \in C^1(\Omega, \R)$ and $t_0 \in I$, then there is a neighborhood $I_0$ around $t_0$ such that
$$\|  \tfrac{d}{dt} F(S(t))  |_{t=t^*} \| \le C \max_{1 \le s \le d; x \in E} \big| \tfrac{ \partial f}{\partial x^s}(x) \big | \qquad \text{ for all } t^* \in I_0,$$
where $C$ is a constant and $E$ is a convex, precompact open set with $\bar{E} \subset \Omega.$ \el

\bp Let $t_0 \in I$ and fix a bounded interval $I_0$ around $t_0$ with $\bar{I_0} \subset I.$ By Theorem \ref{thm2.2}, the joint eigenvalues of $S(t^*)$ are continuous on $I_0$. Thus, there exists an open, precompact, convex set $E \subset \R^d$ such that $ \bar{E} \subset \Omega$ and $\sigma(S(t^*)) \subset E$ for each $t^* \in I_0$. Fix $t^* \in I_0.$ By Proposition \ref{prop3.2}, 
\begin{eqnarray} \label{eqn3.5} \tfrac{d}{dt} F(S(t)) |_{t=t*} = \dis U \Big( \sum_{r=1}^d \ \tilde{\Gamma}^r \tfrac{ \partial f}{\partial x^r}(D) + \left[ Y, f(D) \right] \Big) U^*, \end{eqnarray}
where $U,$ $D^r,$ $\tilde{\Gamma}^r$, and $Y$ are functions of $t^*$ defined in Proposition \ref{prop3.2}, and the joint eigenvalues of $S(t^*)$ are denoted by $x_i$, for $ 1 \le i \le n.$ Observe that the matrix in (\ref{eqn3.5}) can be rewritten as
\begin{eqnarray} \label{eqn3.5.2} \bigg[ \sum_{r=1}^d \ \tilde{\Gamma}^r \tfrac{ \partial f}{\partial x^r}(D) + \left[ Y, F(D) \right] \bigg]_{ij} &=& \left\{
\begin{array}{lr}
\dis \sum_{r=1}^d \Gamma^r_{ij} \tfrac{ \partial f}{\partial x^r}(x_i) & \text{ if } x_i = x_j \\
& \\
\Gamma^q_{ij} \frac{f(x_i)-f(x_j)}{x^q_i -x^q_j} & \text{if\ } x_i \ne x_j, \\ 
\end{array} \right. \end{eqnarray}
where $q$ is such that $x^q_i \ne x^q_j.$ As shown in the proof of Theorem \ref{thm2.1}, the value $ \tfrac{ \Gamma^q_{ij} }{x^q_i - x^q_j}$ is independent of  $q$ whenever  $x^q_i  \ne x^q_j.$  \\

Recall that for a given $n \times n$ self-adjoint matrix $A$ and an $n \times n$ unitary matrix $U$, 
\begin{eqnarray} \label{eqn3.5.1} \max_{ij} | (U  \ A \ U^*) _{ij} |  \le n \|U \  A \ U^*\| =  n\|A \| \le   n^2 \max_{ij} |A_{ij} |.\end{eqnarray}
It is immediate from (\ref{eqn3.5}), (\ref{eqn3.5.2}) and (\ref{eqn3.5.1}) that
\begin{eqnarray} \label{eq3.6} \big | \big | \tfrac{d}{dt} F(S(t)) \big |_{t=t^*} \big | \big | \le n \max \bigg | \sum_{r=1}^d \Gamma^r_{ij} \tfrac{ \partial f}{\partial x^r}(x_i) \bigg | + n \max \bigg | \Gamma^q_{ij} \frac{f(x_i)-f(x_j)}{x^q_i -x^q_j} \bigg |,\end{eqnarray}
where the first maximum is taken over $(i,j)$ with $x_i=x_j$, the second maximum is taken over $(i,j)$ with $x_i \ne x_j,$ and  $q$  is such that $x^q_i \ne x^q_j.$  Fix $(i,j)$ with $x_i \ne x_j.$ Since $f \in C^1(E),$ we can apply the multivariate Mean Value Theorem as follows:
\begin{eqnarray} 
\big | f(x_i) - f(x_j) \big | & = & \big | \nabla f(x^*) \cdot \left( x_i -x_j \right) \big | \nonumber \\ 
& \le & \max_{s;x \in E} \big| \tfrac{ \partial f}{\partial x^s}(x) \big | \dis \sum_{r=1}^d \big |x^r_i -x^r_j \big | \label{eq3.7}, 
\end{eqnarray} 
where $x^*$ is on the line in $E$ connecting $x_i$ and $x_j$. If $x^q_i \ne x^q_j$, for each $r$ with $x_i^r \ne x^r_j$, 
$$\Gamma^q_{ij} \frac{x^r_i -x^r_j}{x^q_i -x^q_j} = \Gamma^r_{ij} .$$ 
It follows from (\ref{eq3.7}) that, for each $(i,j,q)$ with $x^q_i \ne x^q_j$, 
\begin{eqnarray}
\bigg | \Gamma^q_{ij} \frac{f(x_i)-f(x_j)}{x^q_i-x^q_j} \bigg | & \le & \bigg| \tfrac{\Gamma^q_{ij} }{x^q_i -x^q_j} \bigg | \max_{s; x \in E} \big| \tfrac{ \partial f}{\partial x^s}(x) \big | \dis \sum_{r=1}^d \big |x^r_i -x^r_j \big | \nonumber \\
& \le & \ \max_{s; x \in E} \big| \tfrac{ \partial f}{\partial x^s}(x) \big|  \dis \sum_{r=1}^d \big| \Gamma^r_{ij} \big | \nonumber \\
& \le & dn^2 \max_{ s; x \in E} \big| \tfrac{ \partial f}{\partial x^s}(x) \big| \max_{i,j, r} \big | (S^r)'(t^*)_{ij} \big |. \label{eq3.8} 
\end{eqnarray}
Likewise, 
\begin{eqnarray}
\big | \sum_{r=1}^d \Gamma^r_{ij} \tfrac{ \partial f}{\partial x^r}(x_i) \big | & \le & dn^2 \max_{ s; x \in E} \big| \tfrac{ \partial f}{\partial x^s}(x) \big | \max_{i,j,r} \big | (S^r)'(t^*)_{ij} \big|. \label{eq3.9} \end{eqnarray} 
Let $M$ be a constant bounding each $|(S^r)'(t^*)_{ij}|$ on $\bar{I}_0$ and let $C=2dn^3M.$ Substituting (\ref{eq3.8}) and (\ref{eq3.9}) into (\ref{eq3.6}) gives
\beq \hspace{.5in} \big |\big | \tfrac{d}{dt} F(S(t)) \big |_{t=t^*} \big |\big | &\le& 2dn^3 \max_{ s; x \in E} \big| \tfrac{ \partial f}{\partial x^s}(x) \big| \max_{i,j,r} \big | (S^r)'(t^*)_{ij} \big | \\
 & \le&  C \max_{ s; x \in E} \big| \tfrac{ \partial f}{\partial x^s}(x) \big | \qquad  \hspace{1in}\forall \ t^* \in I_0.  \hspace{.2in} \text{\hfill {$\Box$}}\eeq 

\vspace{.1in}

Proof of Theorem \ref{thm3.2}: \\
First assume $\Omega$ is convex. Let $t_0 \in I.$ Let $I_0$ be the interval around $t_0$ and $E$ be the convex, precompact open set given in Lemma \ref{lem3.4}. Since $f$ is a $C^1$ function and $\bar{E}$ is compact, a generalization of the Stone-Weierstrass Theorem in \cite[pg 55]{heb89} guarantees a sequence $ \left \{ \phi_k \right \}$ of functions analytic on $\R^d$ such that
$$ | \phi_k(x) - f(x) | < \tfrac{1}{k} \ \text{ and } \ \big | \tfrac{ \partial \phi_k}{ \partial x^r}(x) - \tfrac{ \partial f}{ \partial x^r}(x) \big | < \tfrac{1}{k} \quad \forall \ x \in \bar{E} \ \text{ and } \ 1 \le r \le d. $$
Lemma \ref{lem3.4} guarantees that, for each $t^* \in I_0,$
\beq \big | \big |  \tfrac{d}{dt} \Phi_k(S(t)) \big |_{t=t^*} - \tfrac{d}{dt} F(S(t)) \big |_{t=t^*} \big | \big | &=& \big | \big | \tfrac{d}{dt}( F-\Phi_k) (S(t)) \big |_{t=t^*}  \big | \big | \\
&& \\
&\le& C \max_{s; x \in E} \big| \tfrac{ \partial (f-\phi_k)}{\partial x^s}(x) \big | \\
&\le& \tfrac{C}{k}, \eeq
where $C$ is a fixed constant. This implies
$$ \left \{ \tfrac{d}{dt} \Phi_k(S(t)) \big |_{t=t^*} \right \} \text{ converges uniformly to } \tfrac{d}{dt} F(S(t)) \big |_{t=t^*} \text{ on } I_0. $$ 
By Proposition \ref{prop3.1}, each $\frac{d}{dt} \Phi_k(S(t)) |_{t=t^*}$ is continuous on $I$. Since the uniform limit of continuous functions is continuous, $\frac{d}{dt} F(S(t)) |_{t=t^*}$ is continuous on $I_0$. \\

Now, let $\Omega$ be an arbitrary domain. Fix $t_0 \in I$ and let $I_0$ be a bounded open interval of $t_0$ with $\bar{I_0} \subset I.$ Let $E \subset \R^d$ be an open precompact set such that $ \bar{E} \subset \Omega$ and $\sigma(S(t^*)) \subset E$ for all $t^* \in I_0.$ Let $O$ be an open set and $K$ be a compact set such that $\bar{E}\subset O \subset K \subset \Omega$ and define a $C^{\infty}$ bump function $b(x)$ such that
\beq
b(x) 
:= \left\{
\begin{array}{ll}
1 & \text{ if } x \in E \\ 
0 & \text{ if } x \in K^c. 
\end{array} \right. 
\eeq
Now we can define a function $g$ in  $C^1(\R^d, \R)$ by
$$g(x):= 
\left\{
\begin{array}{ll}
b(x)f(x) & \text{ if } x \in \Omega \\ 
0 & \text{ if } x \in \Omega^c. 
\end{array} \right. $$
As $\R^d$ is convex, it follows from the previous result that $\frac{d}{dt} G(S(t)) |_{t=t^*}$ is continuous on $I_0$. Since $f(x) = g(x)$ in $E$, It follows from the formula in Proposition \ref{prop3.2} that
$$\tfrac{d}{dt} F(S(t)) |_{t=t^*}= \tfrac{d}{dt} G(S(t)) |_{t=t^*} \qquad \forall \ t^* \in I_0,$$ 
and thus, is continuous in $I_0.$ \ep
\vspace{.1in}

Recall that $CS^d_n$ possesses a Whitney stratification with pieces $\{M_{\alpha}\}.$  Let $\Omega$ be an open set in $\R^d$ and let $f \in C^1(\Omega, \R)$. Let $V$ be an open set in $CS^d_n$ such that  for all $S \in V,$ $\sigma(S) \subset \Omega.$ Define $TV := \cup T(M_{\alpha} \cap V).$ Then, $F(S)$ exists for all $S \in V$  and we can use the derivative results to define a map $DF: TV \rightarrow TS_n.$ \\ 

Specifically, fix an element in $TV,$ which will consist of an $S \in V$ and  $\De \in T_{S}M_{\alpha}$, where $M_{\alpha}$ is the piece containing $S$.  Let $S(t)$ be a $C^1$ curve in $CS^d_n$ such that $S(0)=S$ and $S'(0)= \De.$ Define 
$$DF(S, \De) := \tfrac{d}{dt}F(S(t)) |_{t=0} = \dis U \big ( \sum_{r=1}^d \ \tilde{\Gamma}^r \tfrac{ \partial f}{\partial x^r}(D) + \left[ Y, f(D) \right] \big ) U^*,$$ 
where $U,\ D, \ \tilde{\Gamma}^r,$ and $Y$ are defined using $S$ and $\De$ as in Proposition \ref{prop3.2}. It is easy to see that the map is well-defined and  $DF(S, \cdot)$ is linear in $\De.$  In the following theorem, let $S$ be in a piece $M_{\alpha}$ and let $R$ be in a piece $M_{\beta}$ of a Whitney stratification of $CS^d_n.$

\bt \label{thm3.3} Let $\Omega$ be an open set in $\R^d$ and $V$ be an open set in $CS^d_n$ with $\sigma(S)$ in $\Omega$ for all $S \in V.$ If $f \in C^1(\Omega, \R)$, then 
$$ DF : TV \rightarrow TS_n \text{ is continuous}. $$ 
Specifically, if $S \in V$ with $\Delta \in T_SM_{\alpha}$, then given $\epsilon >0$, there exist $\delta_1, \ \delta_2 >0$ such that if
$R \in V$ with $\Lambda \in T_RM_{\beta}$,  $\| S-R \| < \delta_1,$ and $\| \De - \Lambda \| < \delta_2,$ then
$$\| DF(S,\Delta) -DF(R,\Lambda) \| <  \epsilon.$$ \et

\bp The result for analytic functions follows from Equation (\ref{eq3.1}). For an arbitrary function $f,$ and for $R$ and $\Lambda$ sufficiently close to $S$ and $\Delta$, bound $\| DF(R, \Lambda)\| $ in a manner similar to Lemma \ref{lem3.4}. The remainder of the proof is almost identical to that of Theorem \ref{thm3.2} and is left as an exercise. \ep

\section{ Higher Order Derivatives}
We now consider higher-order differentiation and for ease of notation, discuss only two-variable functions.  We first clarify some  notation. In earlier sections, $(\zeta^1, \dots, \zeta^d)$ referred to a point in $\C^d.$ In this section, $(\zeta_1, \zeta_2)$ denotes a point in $\C^2.$ Previously, $S(t)$ and $T(t)$ denoted two separate curves in $CS^d_n.$ Now, $S(t)$ and $T(t)$ denote the two components of a single curve in $CS^2_n.$ \\ 

Let $(S(t),T(t))$ be a $C^m$ curve in $CS^2_n$ defined on an interval $I.$ If $m \ge 1$, the curve is Lipschitz. By Theorem \ref{thm2.2}, there are Lipschitz curves
\begin{eqnarray} (x_s(t) , y_s(t)) \qquad \text{ for } 1 \le s \le n, \label{eq4.1.2} \end{eqnarray}
defined on $I$  representing the joint eigenvalues of $(S(t),T(t))$. Let $U(t)$ be a unitary matrix diagonalizing $(S(t),T(t))$ so that the joint eigenvalues are ordered as in (\ref{eq4.1.2}).
To simplify notation, we write $(S(t),T(t))$ as $(S,T).$ For $l \in \NN$ with $1 \le l \le m$, define
\begin{eqnarray*} S^l := S^{(l)}(t) \ \ \text{ and }  \ \ T^l := T^{(l)}(t) \label{eq4.1.1a} \end{eqnarray*}
and the set of pairs of index tuples
\beq I_l & := & \left \{ (i_1, \dots, i_k) \cup ( i_{k+1}, \dots, i_j) : i_1 + \dots + i_j =l, i_q \in \NN \text{ for } 1 \le q \le j \right \}. \eeq
For example, $I_2 = \{ (2)\cup \emptyset, \ (1,1) \cup \emptyset, \ (1) \cup (1), \ \emptyset \cup (1,1), \ \emptyset \cup  (2) \}.$ For notational ease, define
$$U := U(t), \ \ \  x_s :=x_s(t), \ \ \  y_s :=y_s(t) \ \ \ \text{ for } 1 \le s \le n.$$ 
For some formulas, we will conjugate the derivatives in (\ref{eq4.1.1a}) by $U^*$ and so define
$$ \Gamma^l := U^* \ S^{(l)} \ U \ \ \  \text{ and } \ \ \ \Delta^l := U^* \ T^{(l)} \ U,  \ \ \ \text{ for }  1 \le l \le m.$$ 
We will use the integral formula given in Lemma \ref{lem3.1} and simplify it by defining
$$R_1 := (\zeta_1I - S)^{-1} \ \ \  \text{ and } \ \ \ R_2 :=(\zeta_2I - T)^{-1},$$ 
where $\zeta_1$ and $\zeta_2$ are in the resolvents of $S$ and $T$ respectively. Now, let $J_1$ and $J_2$ be open intervals in $\R$ and let $f$ be an element of $C^m( J_1 \times J_2, \R).$ Fix $j$ and $k$ in $\NN$ such that $k \le j \le m$. Fix $k+1$ points $x_1, \dots, x_{k+1}$ in $J_1$ and $j-k+1$ points $y_1, \dots, y_{j-k+1}$ in $J_2$. Then
$$f^{[k,j-k]}(x_1, \dots, x_{k+1}; y_1, \dots, y_{j-k+1})$$
denotes the divided difference of $f$ taken in the first variable $k$ times and the second variable $j-k$ times, evaluated at the given points.  Finally, let $\odot$ denote the Schur product of two matrices. We will prove the following differentiability result:

\bt \label{thm4.1} Let $J_1$ and $J_2$ be open intervals in $\R$ and let $f \in C^m(J_1 \times J_2, \R).$ Let $(S,T)$ be a $C^m$ curve in $CS^2_n$ defined on an interval $I$ with joint eigenvalues in $J_1 \times J_2.$  For $1 \le l \le m$ and $t^* \in I,$ $\frac{d^l}{dt^l} F(S,T)|_{t=t^*}$ exists and 
\beq \tfrac{d^l}{dt^l} F(S,T) \big|_{t=t^*} = U\bigg( \sum_{I_l} \sum_{s_2. . s_j =1}^{n} \frac{l!}{i_1! \cdots i_j!} \left[ f^{[k,j-k]}(x_{s_1}, \dots, x_{s_{k+1}}; y_{s_{k+1}}, \dots,y_{s_{j+1}}) \right]_{s_1, s_{j+1}=1}^n \\
\  \ \odot \ \left[ \Gamma^{i_1}_{s_1s_2}\dots \Gamma^{i_k}_{s_ks_{k+1}} \Delta^{i_{k+1}}_{s_{k+1}s_{k+2}} \dots \Delta^{i_j}_{s_js_{j+1}} \right]_{s_1, s_{j+1}=1}^n \bigg) U,^* \hspace{1in} \eeq 
where the $U$, $U^*$, $\Gamma^i$, $\Delta^j$, $x_q$ and $y_r$ are evaluated at $t^*.$\et
Notice that the derivative formula in Theorem \ref{thm4.1} requires $f$ to be defined on pairs $(x_q,y_r)$ for $ 1 \le r,q \le n$, rather than just at the joint eigenvalues $(x_q,y_q)$ of $(S,T).$  This condition was not needed in Theorem \ref{thm3.1}.  Before proving Theorem \ref{thm4.1}, we consider the case where $f$ is real-analytic and show:

\bprop \label{prop4.1} Let $J_1$ and $J_2$ be open intervals in $\R$ and let $f$ be real-analytic on $J_1 \times J_2$. Fix $m \in \NN$ and let $(S,T)$ be a $C^m$ curve in $CS^2_n$ defined on an interval $I$ with joint eigenvalues in $J_1 \times J_2.$ Then
$\tfrac{d^m}{dt^m} F(S,T)$ exists, has the form in Theorem \ref{thm4.1}, and $\tfrac{d^m}{dt^m} F(S,T)  |_{t=t^*} \text{ is continuous as a function of } t^* \text{ on } I .$ \eprop

The proof of Proposition \ref{prop4.1} requires the following two technical lemmas:

\bl \label{lem4.1} Let $(S,T)$ be a $C^m$ curve in $CS^2_n$ defined on an interval $I$. Let $t^* \in I$ and let $\zeta_1$ and $\zeta_2$ be elements in the resolvents of $S(t^*)$ and $T(t^*)$ respectively. Then 
\beq  \tfrac{d^{l}}{dt^{l}} \big ( R_1R_2 \big) \big |_{t=t^*}= \sum_{I_l} \frac{l!}{i_1! \cdot \cdot \ i_j!} R_1S^{i_1}R_1 \dots S^{i_k}R_1R_2T^{i_{k+1}}R_2 \dots T^{i_j}R_2  \  \ \forall \ 1\le l \le m,\eeq
where each $R_1,$  $R_2,$ $S^i,$ and $T^j$ is evaluated at $t^*.$
\el
\bp The proof is a technical calculation using induction on $l$ and the formulas 
$$ \hspace{1.5in} \tfrac{d}{dt}R_1 = R_1S^1R_1 \ \ \text{ and } \ \ \tfrac{d}{dt}R_2 = R_2T^1R_2 . \hspace{1.35in} \text{\hfill{$\Box$}} $$

\bl \label{lem4.2} Let $J_1$ and $J_2$ be open intervals in $\R$ and let $f$ be real-analytic on $J_1 \times J_2$. Let $j \ge k \in \NN$. Choose $k+1$ points $x_1, \dots, x_{k+1} \in J_1$ and $j-k+1$ points $y_1, \dots, y_{j-k+1} \in J_2.$ Extend $f$ to be analytic on a complex rectangle $\tilde{\Omega} \subset \C^2$ such that each $(x_q,y_r) \in \tilde{\Omega}.$ Then $f^{[k, j-k]}(x_1, \dots, x_{k+1}; y_1, \dots, y_{j-k+1})$ exists and
\beq f^{[k, j-k]}(x_1,.., x_{k+1}; y_1,.., y_{j-k+1})   = \frac{1}{(2\pi i)^2}\dis \int_{C_2} \int_{C_1} \dfrac{f(\zeta_1, \zeta_2)} {\prod_{q=1}^{k+1} (  \zeta_1-x_q)  \prod_{r=1}^{j-k+1}(\zeta_2-y_r)}  d\zeta_1 d\zeta_2,\eeq
where $C_1$ and $C_2$ are rectifiable curves strictly enclosing $x_1, \dots, x_{k+1}$ and $y_1, \dots, y_{j-k+1}$ respectively, such that $C_1 \times C_2 \subset \tilde{\Omega}.$ \el
\bp  For a one-variable function, the formula is proven in \cite[pg 2]{don74} and the  two-variable analogue follows easily from the one variable case. \ep

\vspace{.1in}

Proof of Proposition \ref{prop4.1}: \\
Use the integral formula in Lemma \ref{lem3.1} to establish an integral formula for $ \tfrac{d^m}{dt^m} F(S,T)$ similar to the first line of (\ref{eq3.1}). Simplify the formula using Lemma \ref{lem4.1}. This formula implies that the derivative is continuous.  Then, use Lemma \ref{lem4.2} to convert the derivative into a formula involving the divided differences of $f$. The details are left as an exercise. \ep 

\vspace{.1in}

Proof of Theorem \ref{thm4.1}: \\
The result follows via induction on $l$, and the base case is covered by Theorem \ref{thm3.1}. For the inductive step, fix $t^* \in I$. Let $p$ be a polynomial such that $p$ and its derivatives to $l^{th}$ order agree with $f$ at the points $(x_{q}(t^*),y_r(t^*))$ for $1 \le q,r \le n.$ Find a constant $C$ such that for $t$ near $t^*$,
$$ \| \tfrac{d^{l-1}}{dt^{l-1}} F(S,T)  -\tfrac{d^{l-1}}{dt^{l-1}} P(S,T)  \|  \le C \max  \big |  (f-p)^{[k,j-k]}(x_{s_1}, .. , x_{s_{k+1}}; y_{s_{k+1}}, .. , y_{s_{j+1}}) \big|, $$
where the joint eigenvalues of $(S,T)$ are given by $(x_q,y_q)$ and the maximum is over  $(k,j)$ with $ k \le j <l\in \NN$ and sets $\left \{ (s_1, . . , s_{k+1}) \cup (s_{k+1}, .. , s_{j+1}) : 1 \le s_{1}..s_{j+1}\le n \right \}.$  As in Theorem \ref{thm3.1}, apply the multivariate Mean Value theorem to each $(f-p)^{[k,j-k]}$ and  use the Lipschitz property of the eigenvalues to conclude
$$\tfrac{d^l}{dt^l} F(S,T) |_{t=t^*} \text{ exists and equals } \tfrac{d^l}{dt^l} P(S,T) |_{t=t^*}.$$
The details are left as an exercise. \ep

\vspace{.1in}

We now show that the formula in Theorem \ref{thm4.1} is continuous.

\bt \label{thm4.2} Let $J_1$ and $J_2$ be open intervals in $\R$ and $f \in C^m(J_1 \times J_2, \R)$. Let $(S,T)$ be a $C^m$ curve in $CS^2_n$ defined on an interval $I$ with joint eigenvalues in $J_1 \times J_2.$ Then for all $l \in \NN$ with $1 \le l \le m$, 
$$\tfrac{d^l}{dt^l} F(S,T) |_{t=t^*} \text{ is continuous as a function of } t^* \text{ on } I.$$
\et

For the proof, we require the following lemma. The result is well-known for one-variable functions, and Brown and Vasudeva prove this two-variable analogue in \cite{brv00}: 

\bl \label{lem4.5.1} Let $J_1$ and $J_2$ be open intervals in $\R$ and let $f \in C^m(J_1 \times J_2, \R).$ Choose $j, k \in \NN$ with $ k \le j \le m.$ Let $x_1, \dots, x_{k+1} \in J_1$ and $y_1, \dots, y_{j-k+1} \in J_2$ and choose closed subintervals $\tilde{J_1}$ and $\tilde{J_2}$ containing the $x$ and $y$ points respectively. Then, there exists $(x^*, y^*) \in \tilde{J_1} \times \tilde{J_2}$ with
$$f^{[k,j-k]}(x_1, \dots, x_{k+1}; y_1, \dots, y_{j-k+1}) = \frac{f^{(k,j-k)}(x^*, y^*)}{k!(j-k)!}.$$ \el

Proof of Theorem \ref{thm4.2}:\\
For $l <m$, the result follows from Theorem \ref{thm4.1}, which implies that $\frac{d^l}{dt^l} F(S,T)$ is differentiable and, hence, continuous. \\

For $l =m$, fix $t_0 \in I.$  As in Lemma \ref{lem3.4}, find a constant $C$ and an open precompact convex set $J$ with $\bar{J} \subset J_1 \times J_2$ such that, for all $g \in C^m(J_1 \times J_2, \R)$ and $t^*$ near $t_0,$ 
\beq \big | \big |  \tfrac{d^m}{dt^m} G(S,T) |_{t=t^*} \big | \big| \le  C  \ \max_{ \left \{ j, k \ ;\ (x,y)  \in J \right \}} |f^{(k,j-k)}(x,y)|, \eeq
where $0 \le k \le j \le m.$ The estimates for this bound require Lemma \ref{lem4.5.1}. Then, approximate $f$ to $m^{th}$ order uniformly on $\bar{J}$  by analytic functions $\{\phi_r \}$ and show
$$\{ \tfrac{d^m}{dt^m} \Phi_r (S,T)|_{t=t^*} \} \text{ converges uniformly to }  \tfrac{d^m}{dt^m} F (S,T)|_{t=t^*}$$
in a neighborhood of $t_0$. The result then follows from Proposition \ref{prop4.1}. \ep

\section{ Applications }
The formulas in Proposition \ref{prop3.2} and Theorem \ref{thm4.1} can be used to analyze monotonicity and convexity of matrix functions.  A function $F: S_n \rightarrow S_n$  is \textbf{matrix monotone} if 
\begin{eqnarray*} F(A) \ge F(B) \text{ whenever } A \ge B \qquad \forall \ A, B \in S_n. \end{eqnarray*} 
For $F$ continuously differentiable, an equivalent condition is 
\begin{eqnarray} \tfrac{d}{dt}F(S(t))|_{t=t^*} \ge 0 \text{ whenever } S'(t^*) \ge 0, \qquad \forall \ C^1 \ S(t) \in S_n. \label{eqn5.1} \end{eqnarray}
The local monotonicity condition in (\ref{eqn5.1}) extends to multivariate matrix functions: the only adjustment is that $S(t)$ is in $CS^d_n$. In \cite{mcc10}, Agler, McCarthy, and Young characterized such  locally monotone matrix functions on $CS^d_n$ using a special case of Theorem \ref{thm3.1} and Proposition \ref{prop3.2}. Specifically, they had to assume that $S(t)$ had distinct joint eigenvalues at each $t$. Our results in Section 3 extend the derivative formula to general $C^1$ curves in $CS^d_n$ and show that the formula is continuous. \\

A matrix function $F:S_n \rightarrow S_n$ is \textbf{matrix convex} if 
\begin{eqnarray} F( \lambda A + (1-\lambda) B) \le \lambda F(A) + (1-\lambda) F(B) \ \ \ \forall \ A,B \in  S_n, \ \lambda \in [0,1] \label{eqn5.2}. \end{eqnarray}
This condition extends to multivariate matrix functions with an additional restriction on the pairs $A,B$ in $CS^d_n$; we also require $\lambda A + (1-\lambda) B \in CS^d_n$  for $\lambda \in [0,1].$  Given such $A,B$, define $S(t)$ on $[0,1]$ by
$$S^r(t):=  t A^r + (1-t) B^r \qquad \forall \ 1 \le r \le d.$$
If $F$ is twice continuously differentiable, it can be shown that (\ref{eqn5.2}) is equivalent to 
\begin{eqnarray}  \tfrac{d^2}{dt^2} F(S(t))|_{t=t^*} \ge 0 \qquad \forall \text{ such } S(t), \ t^* \in [0,1]. \label{eqn5.3} \end{eqnarray} 
Assume $F$ was defined using a real function $f$ as in (\ref{eqn1.1}). For $d=2,$ Theorem \ref{thm4.1} tells us that, up to conjugation by a unitary matrix $U$, 
\begin{eqnarray*} \Big[ \tfrac{d^2}{dt^2} F(S(t))|_{t=t^*} \Big]_{ij} &=&  2 \sum_{k=1}^n f^{[2,0]}(x_i,x_k,x_j;y_j) \Gamma_{ik}\Gamma_{kj} + f^{[1,1]}(x_i,x_k;y_k,y_j) \Gamma_{ik}\Delta_{kj} \nonumber \\
& & \hspace{.25in} + \  f^{[0,2]}(x_i;y_i,y_k,y_j) \Delta_{ik}\Delta_{kj}, \label{eqn5.4} \end{eqnarray*}
where $\{(x_i,y_i) : 1 \le i \le n \}$ are the joint eigenvalues of  $t^*A + (1-t^*)B$ and 
$$\Gamma: = U^* \big(   A_1-B_1 \big) U \ \text{ and } \Delta: = U^* \big( A_2-B_2 \big) U.$$ 
This formula can be simplified using the relationship between $\Gamma$ and $\Delta$  discussed in Theorem \ref{thm2.1}.  Specifically, we know
$$  (x_i - x_j)\Delta_{ij} = (y_i - y_j)\Gamma_{ij}  \qquad \forall \ 1 \le i,j \le n. $$  
Thus, this formula gives a characterization of convex matrix functions on $CS^2_n$.

\bibliography{references2}
\end{document}